\documentclass[a4paper,11pt,reqno]{amsart}
\usepackage{amsfonts}
\usepackage{amssymb}
\usepackage{mathrsfs}
\usepackage{a4wide}
\usepackage{hyperref}
\usepackage{theoremref}
\usepackage[dvipsnames]{xcolor}
\usepackage[utf8]{inputenc}
\usepackage{array}
\usepackage{enumitem}
\usepackage{comment}
\usepackage{graphicx,tikz,parskip}
\usepackage{csquotes}
\usepackage{subfiles}
\usepackage{appendix}
\usepackage{rotating}
\usepackage{xcolor}
\usepackage{float}
\usepackage{mathtools}
\usepackage{pifont}
\usepackage{textcomp}
\usepackage{multirow}
\usepackage{cancel}
\usepackage{todonotes}

\expandafter\def\expandafter\normalsize\expandafter{%
    \normalsize%
    \setlength\abovedisplayskip{3pt}%
    \setlength\belowdisplayskip{3pt}%
    \setlength\abovedisplayshortskip{0pt}%
    \setlength\belowdisplayshortskip{0pt}%
}

\usepackage{xpatch}
\makeatletter
\xpatchcmd{\proof}{\topsep6\p@\@plus6\p@\relax}{}{}{}
\makeatother

\theoremstyle{plain}
\newtheorem{lemma}{Lemma}[section]
\newtheorem{theorem}[lemma]{Theorem}
\newtheorem{corollary}[lemma]{Corollary}
\newtheorem{proposition}[lemma]{Proposition}
\newtheorem*{theorem*}{Theorem}
\newtheorem{counterexample}[lemma]{Counterexample}

\theoremstyle{definition}
\newtheorem{definition}[lemma]{Definition}
\theoremstyle{remark}
\newtheorem{remark}[lemma]{Remark}

\newtheorem{example}[lemma]{Example}

\newtheorem{manualremark}{Remark}
\newenvironment{remark*}[1]{%
  \renewcommand\themanualremark{#1}%
  \manualremark
}{\endmanualremark}

\renewcommand{\l}{\langle}
\renewcommand{\r}{\rangle}
\renewcommand{\phi}{\varphi}
\renewcommand{\bf}{\textbf}

\makeatletter
\newcommand{\typeitem}[1]{\item[#1]\protected@edef\@currentlabel{#1}}
\makeatother

\makeatletter
\@namedef{subjclassname@2020}{\textup{2020} Mathematics Subject Classification}
\makeatother

\begin{document}

\title[Weakly right coherent monoids]{Weakly right coherent monoids}

\author{Levent Michael Dasar, Victoria Gould, Craig Miller}
\address{Department of Mathematics, University of York, York YO10 5DD, UK}
\email{levent.dasar@york.ac.uk; victoria.gould@york.ac.uk; craig.miller@york.ac.uk}



\begin{abstract}
A monoid $S$ is said to be weakly right coherent if every finitely generated right ideal of $S$ is finitely presented as a right $S$-act. It is known that $S$ is weakly right coherent if and only if it satisfies the following conditions: $S$ is right ideal Howson, meaning that the intersection of any two finitely generated right ideals of $S$ is finitely generated; and the right annihilator congruences of $S$ are finitely generated as right congruences.  We examine the behaviour of these two conditions (in the more general setting of semigroups) under certain algebraic constructions and deduce closure results for the class of weakly right coherent monoids. In particular, we show that the  class of weakly right coherent monoids is closed under taking direct products and monoid free products. Finally, we provide a survey of connections between each of the aforementioned conditions and axiomatisability of classes of acts.
\end{abstract}

\maketitle

 \section{Introduction}
An algebra from a given variety is said to be {\em coherent} if all its finitely generated subalgebras are finitely presented.  In the context of groups and rings, coherency has received significant attention. For a survey of coherency in groups, see \cite{wise}.  A ring is called {\em right coherent} if all its finitely presented right modules are coherent.  It turns out that right coherency for rings is equivalent to every finitely generated right ideal being finitely presented as a right module \cite[Theorem 2.1]{chase:1960}. 
For monoids, acts play the analagous role to that of modules in the theory of rings.  Thus, a monoid $S$ is called {\em right coherent} if every finitely presented right $S$-act is coherent.  Right coherency for monoids was introduced in \cite{coherentmonoids} and has since been intensively studied; see \cite{bgr:2023,d:2020,dg:2023,Hartmann,freemonoid}.  In particular, it has been shown that groups, free monoids and the bicyclic monoid are right (and left) coherent. 

A monoid is {\em weakly right coherent} if all its finitely generated right ideals are finitely presented as right acts.  In contrast to the situation for rings, right coherency for monoids is strictly stronger than weak right coherency.  For instance, infinite full transformation monoids are weakly right (and left) coherent \cite[Corollary 3.6, Theorem 3.7]{coherentmonoids} but not right (or left) coherent \cite[Theorem 3.4]{bgr:2023}.  

It was shown in \cite[Corollary 3.3]{coherentmonoids} that a monoid $S$ is weakly right coherent if and only if the following conditions hold: the intersection of any two finitely generated right ideals of $S$ is finitely generated; and the right annihilator congruences $\bf{r}_S(a)=\{(s,t)\in S\times S \mid as=at\}$ ($a\in S$) of $S$ are finitely generated.  These two conditions can of course be defined, more generally, for semigroups.  We call semigroups satisfying the latter condition {\em finitely right equated}.  
Semigroups satisfying the former property are known as {\em right ideal Howson} and were studied in \cite{Carson}.  We remark that a monoid is right ideal Howson if and only if it is finitely aligned \cite{exel:18}; for semigroups, being finitely aligned is a stronger condition, as  demonstrated in \cite{Carson}.

Right coherency and weak right coherency for monoids have been formulated in terms of the model-theoretic notion of axiomatisability.  A class $C$ of first-order structures is {\em axiomatisable} if there is a set of sentences whose models are exactly the members of $C.$
A classical result of Wheeler \cite{wheeler:1976}, interpreted for acts in \cite[Theorem 6]{modelcompanions}, says that a monoid is right coherent if and only if the class of existentially closed right acts is axiomatisable.  By \cite[Theorem 3, Corollary 4.2]{modelcompanions}, a monoid is weakly right coherent if and only if the class of $\aleph_0$-injective right acts is axiomatisable.

The main purpose of the present article is to investigate the closure properties of the class of weakly right coherent monoids with respect to some standard algebraic and semigroup-theoretic constructions.  Since some of the constructions considered yield semigroups that are not monoids, it will be useful to have a notion of weak right coherency for semigroups.  Thus, we define a semigroup $S$ to be {\em weakly right coherent} if $S^1$ is weakly right coherent.  Our strategy is to first study the closure properties of the classes of right ideal Howson semigroups and finitely right equated semigroups, and then deduce results concerning weak right coherency.  In the case of right ideal Howson semigroups, our work builds on that of \cite[Section 4]{Carson}.
Our main results are summarised in Figure \ref{table}.

The paper is structured as follows. In Section \ref{sec:prelim} we introduce the necessary preliminary material on semigroups. Sections \ref{sec:rih} and \ref{sec:fre} are concerned with the properties of being right ideal Howson and being finitely right equated, respectively. In Section \ref{sec:wrc} we use the results of the previous two sections to derive closure results for the class of weakly right coherent monoids.  Finally, in Section \ref{sec:axiom}, we provide a survey of connections between axiomatisability of certain classes of acts and each of the properties of being right ideal Howson, being finitely right equated, and weak right coherency.

\section{Preliminaries}\label{sec:prelim}
In this section we introduce the necessary preliminary material. For a more complete review of basic semigroup theory we refer the reader to \cite{Howie}. 
Throughout this section, $S$ denotes a semigroup.  We denote by $S^1$ ($S^0$)  the monoid (semigroup) obtained from $S$ by adjoining an identity (zero) if necessary. 

Recall that a right ideal $I$ of $S$ is \textit{finitely generated} if there exists a finite set $X\subseteq I$ such that $XS^{1}=I$.
A \textit{right congruence} $\rho$ on $S$ is an equivalence relation on $S$ such that $(a,b)\in\rho$ implies $(as,bs)\in \rho$ for all $s\in S$. For a set $X\subseteq S\times S$, we denote the smallest right congruence on $S$ containing $X$ by $\l X\r$. A right congruence $\rho$ is \textit{finitely generated} if there exists a finite set $X\subseteq\rho$ such that $\rho=\l X\r$. The following lemma, a particular case of a more general fact from universal algebra, provides an explicit characterisation of $\l X\r$.

\begin{lemma}\label{X-sequence}
Let $S$ be a semigroup and $X\subseteq S\times S.$
For $a,b\in S,$ we have $(a,b)\in\l X\r$ if and only if either $a=b$ or there is a sequence
\[a=p_1c_1,\ q_1c_1=p_2c_2,\ \dots\ ,\ q_nc_n=b\]
where $(p_i,q_i)\in X$ or $(q_i,p_i)\in X,$ and $c_i\in S^{1}$, for each $1\leq i\leq n\in\mathbb{N}$.
\end{lemma}\vspace{-2mm}

A sequence of the form given in Lemma \ref{X-sequence} is called an {\em $X$-sequence from $a$ to $b$}.  For convenience, and without loss of generality, from now on we will assume, unless stated otherwise, that a generating set $X$ of a right congruence is symmetric, in the sense that $(p,q)\in X$ if and only if $(q,p)\in X$.  

The \textit{Rees quotient} of $S$ by an ideal $I$ of $S$, denoted $S/I$, is the semigroup with underlying set $S{\setminus}I\cup\{0\},$ where $0\notin S{\setminus}I,$ and multiplication defined by $a\cdot b=ab$ if $a,b,ab\in S{\setminus}I$, or else $a\cdot b=0.$
A \textit{retract} of $S$ is a subsemigroup $T$ of $S$ such that there is a homomorphism $\phi:S\rightarrow T$ with $t\phi =t$ for all $t\in T$.  If a subsemigroup $T$ of $S$ is such that $S{\setminus}T$ is finite, we say that $T$ is \textit{large} in $S$ and that $S$ is a \textit{small extension} of $T.$

Let $\{S_{i}\}_{i\in I}$ be a collection of pairwise disjoint semigroups and $\sigma:\bigsqcup_{i\in I}S_i\rightarrow I$ the source function given by $s\sigma=i$ whenever $s\in S_i$. The \textit{semigroup free product} $\prod^{\ast}\{S_{i}\mid i\in I\}$ of $\{S_{i}\mid i\in I\}$ is the semigroup whose elements are all sequences $s_1\ast\dots\ast s_{n}$ where $s_j\in \bigsqcup_{i\in I}S_i$ for $1\leq j\leq n$ and $s_k\sigma\neq s_{k+1}\sigma$ for $1\leq k< n$, with multiplication defined by
\[
(s_1\ast\dots\ast s_n)\ast(t_1\ast\dots\ast t_m)=\begin{cases}
    s_1\ast\dots\ast s_n\ast t_1\ast\dots\ast t_m & \text{if $s_n\sigma\neq t_1\sigma$,}\\
    s_1\ast\dots\ast s_nt_1\ast\dots\ast t_m & \text{otherwise.}
\end{cases}
\]
In view of the associativity of the operation,  we  consider {\em any} sequence $s_1\ast\dots\ast s_{n}$, where $s_j\in \bigsqcup_{i\in I}S_i$ for $1\leq j\leq n$,
as representing an element of the free product,  noting that there is a unique representation of an element of the free product in this form where additionally $s_k\sigma\neq s_{k+1}\sigma$ for $1\leq k< n$; such a representation is called {\em reduced}. 

Suppose that each semigroup $S_i$ above is a {\em monoid} with identity $1_i$. Let $\rho$ be the congruence on $\prod^{\ast}\{S_{i}\mid i\in I\}$ generated by $\{(1_i,1_j)\mid i,j\in I\}$. The \textit{monoid free product} $\prod_1^{\ast}\{S_i\mid i\in I\}$ of $\{S_i\mid i\in I\}$ is $\prod^{\ast}\{S_i\mid i\in I\}/\rho$.  For convenience, as in the semigroup case, we denote elements of $\prod_1^{\ast}\{S_i\mid i\in I\}$ by sequences 
$s_1\ast\dots\ast s_{n}$ where $s_j\in \bigsqcup_{i\in I}S_i$ for $1\leq j\leq n$,   noting again that such expressions are not unique. In each non-identity $\rho$-class there exists a unique sequence such that $s_k\sigma\neq s_{k+1}\sigma$ for $1\leq k< n$ and which contains no entries from $\{1_i\mid i\in I\}$; we again call such a sequence \textit{reduced}.

\begin{figure}[h]
\begin{tabular}{|c|ccl|}
\hline
\multirow{2}{*}{\bf{Construct}} & \multicolumn{3}{c|}{\bf{Preservation}}                                                                                                                      \\ \cline{2-4} 
& \multicolumn{1}{c|}{RIH} & \multicolumn{1}{c|}{\quad FRE\quad\quad}                            & \multicolumn{1}{c|}{\quad WRC\quad\quad}               \\ \hline\hline
Quotient                            & \multicolumn{1}{l|}{\ding{55} \,\cite[Prop.\! 4.1]{Carson}}  & \multicolumn{1}{l|}{\ding{55} \,[Cex.\! \ref{fre:rq}]}                                 & \ding{55} \,[Cex.\! \ref{wrc:rq}]                              \\ \hline
Rees quotient                       & \multicolumn{1}{l|}{\ding{51} [Cor.\! \ref{rih:rq}]}                                  & \multicolumn{1}{l|}{\ding{55} \,[Cex.\! \ref{fre:rq}]}                                 & \ding{55} \,[Cex.\! \ref{wrc:rq}]                                 \\ \hline
Retract                 & \multicolumn{1}{l|}{\ding{51} [Cor.\! \ref{retract}]}   & \multicolumn{1}{l|}{\ding{51} [Cor.\! \ref{retractR}]}                                 &       \ding{51} [Thm.\! \ref{wrc:retract}]                            \\ \hline
(Large) subsemigroup                           & \multicolumn{1}{l|}{\ding{55} \,[Cex.\! \ref{rih:lsse}]}                                 & \multicolumn{1}{l|}{\ding{55} \,[Cex. \!\ref{fre:lsse}]}                                 & \ding{55} \,[Cex.\! \ref{wrc:lsse}]                                 \\ \hline
Ideal                 & \multicolumn{1}{l|}{\ding{55} \,[Cex.\! \ref{Howson,ideal}]}   & \multicolumn{1}{l|}{\ding{55} \,[Cex.\! \ref{s1notsR}]}                                  &    \ding{55} \,[Rem.\! \ref{wrc:ideal}]                               \\ \hline
Small extension                           & \multicolumn{1}{l|}{\ding{55} \,[Cex.\! \ref{rih:lsse}]}                                 & \multicolumn{1}{l|}{\ding{55} \,[Cex. \!\ref{fre:lsse}]}                                 & \ding{55} \,[Cex.\! \ref{wrc:lsse}]                                 \\ \hline
Adjoin an identity                           & \multicolumn{1}{l|}{\ding{51} \cite[Lem. 2.2]{Carson}}                                 & \multicolumn{1}{l|}{\ding{51} [Cor. \!\ref{adjoin1R}]}                                 & \ding{51} [Thm.\! \ref{wrc}]                                 \\ \hline
Adjoin a zero                           & \multicolumn{1}{l|}{\ding{51} [Cor.\! \ref{adjoin0}]}                                 & \multicolumn{1}{l|}{\ding{55} \,[Rem. \ref{fre:not0}\!]}                                 & \ding{51} [Thm.\!  \ref{wrc:0}]                                 \\ \hline
Semigroup direct product                     & \multicolumn{1}{l|}{\ding{55} \,\cite[Prop.\! 4.4]{Carson}}   & \multicolumn{1}{l|}{\ding{55} \,[Cex.\! \ref{fre:sgrpdp}]}                                  &              \ding{55} \,[Cex.\! \ref{wrc:sgrpdp}]                 \\ \hline
Monoid direct product                     & \multicolumn{1}{l|}{\ding{51} \cite[Prop.\! 4.5]{Carson}}   & \multicolumn{1}{l|}{\ding{51} [Cor.\! \ref{fre:prod}]}                                  & \ding{51} [Thm.\! \ref{wrc:mdp}]                                 \\ \hline
Semigroup free product                 & \multicolumn{1}{l|}{\ding{51} \cite[Prop.\! 4.2]{Carson}}   & \multicolumn{1}{l|}{\ding{55} \,[Cex.\! \ref{cex:sfp}]}                                  &       \ding{51} [Thm.\! \ref{wrc:sfp}]                         \\ \hline
Monoid free product                 & \multicolumn{1}{l|}{\ding{51} \cite[Prop.\! 4.3]{Carson}}   & \multicolumn{1}{l|}{\ding{51} [Thm.\! \ref{fre:mfp}]}                                  & \ding{51} [Thm.\! \ref{wrc:mfp}]                          \\ \hline
\end{tabular}
\caption{Closure properties of the classes of right ideal Howson (RIH) semigroups, finitely right equated (FRE) semigroups and weakly right coherent (WRC) semigroups, where \ding{51} indicates closure and \ding{55} indicates non-closure of the class under performing the construction.}\label{table}
\end{figure}

\section{Right ideal Howson semigroups}\label{sec:rih}

\begin{definition}
A semigroup is \textit{right ideal Howson} (\textit{RIH}) if the intersection of any two finitely generated right ideals is finitely generated.    
\end{definition}\vspace{-2mm}

\begin{proposition}\cite[Lemma 2.2]{Carson}\label{adjoin1} A semigroup $S$ is RIH if and only if $S^1$ is RIH.
\end{proposition}\vspace{-2mm}

A  semigroup is RIH if and only if the intersection of principal right ideals is finitely generated. Also, since the intersections of principal right ideals generated by $\mathcal{R}$-comparable elements are clearly principal,  we need not consider these. We will use both of these facts without explicit mention throughout this paper.

We begin by considering the closure of the class of RIH semigroups under algebraic constructions, building on the earlier work of Carson and Gould \cite{Carson}.

Since free semigroups are RIH, but not all semigroups are RIH \cite{freemonoid}, certainly the property of being RIH is not closed under homomorphic images. However, this property is closed under Rees quotients and retracts, as a consequence of the following technical result.

\begin{proposition}\label{commutes}
Let $S$ and $T$ be semigroups such that for all $\mathcal{R}$-incomparable elements $u,v\in T$ there exist $a,b\in S$ and a homomorphism $\phi : S\to T$ such that $a\phi=u$, $b\phi=v$ and $(aS\cap bS)\phi=uT\cap vT$. If $S$ is RIH then so is $T$.
\begin{proof}
Let $u,v\in T$ be $\mathcal{R}$-incomparable. Then, by assumption, there exist $a,b\in S$ and a homomorphism $\phi:S\rightarrow T$ such that $a\phi=u$, $b\phi=v$ and $(aS\cap bS)\phi=uT\cap vT$. Then $a,b\in S$ are $\mathcal{R}$-incomparable, and hence $aS^1\cap bS^1=aS\cap bS$.  As $S$ is RIH, we have $aS\cap bS=XS^1$ for some finite set $X$, and by application of $\phi$ we obtain 
$$uT^1\cap vT^1=uT\cap vT=(aS\cap bS)\phi=(XS^1)\phi
\subseteq (X\phi)T^1\subseteq 
(aS\cap bS)\phi T^1\subseteq uT^1\cap vT^1 .$$ Therefore 
$uT^1\cap vT^1=(X\phi)T^1$ and $T$ is RIH.
\end{proof}
\end{proposition}\vspace{-2mm}

\begin{corollary}\label{rih:rq}
The class of RIH semigroups is closed under Rees quotients.
\begin{proof}    
Given a semigroup $S$ and ideal $I$ of $S$, it is straightforward to show that the conditions of Proposition \ref{commutes} are satisfied, where for any $u,v\in S/I$ we use the canonical homomorphism $\phi:S\to S/I.$
\end{proof}
\end{corollary}\vspace{-2mm}

\begin{corollary}\label{retract}The class of RIH semigroups is closed under retracts.
\begin{proof}
Let $S$ and $T$ be semigroups with a retraction $\phi:S\rightarrow T.$  Consider $\mathcal{R}$-incomparable $u,v\in T$. Then $u\phi=u$ and $v\phi=v,$ and certainly $(uS\cap vS)\phi\subseteq uT\cap vT$ (since $\phi$ is a homomorphism).  For $y\in uT\cap vT,$ applying $\phi$ we obtain $y=y\phi\in(uS\cap vS)\phi$.  Thus $(uS\cap vS)\phi=uT\cap vT$. The result now follows from Proposition \ref{commutes}.
\end{proof}
\end{corollary}\vspace{-2mm}

Next, we show that the property of being RIH is not in general inherited by a large subsemigroup or a small extension.

\begin{counterexample}\label{rih:lsse}
The class of RIH semigroups is not closed under small extensions or large subsemigroups.
    \begin{proof}
We first show that the class of RIH semigroups is not closed under small extensions. Let $F$ be the free commutative semigroup generated by $\{x_i\mid i\in\mathbb{N}\}$ and let $\overline{F}=\{\overline{w}\mid w\in F\}$. Then $T=F\cup\overline{F}\cup\{0\}$ forms a semigroup with zero where the operation on $F$ is extended as follows:
\[
u\overline{v}=\overline{u}v=\overline{uv}\quad\text{ and }\quad\overline{u}\,\overline{v}=a0=0a=0\quad(u,v\in F, a\in T).
\]
We claim that $T$ is RIH. Let $u\in F$. Then $uT^1=uF^{1}\cup\{\overline{uw}\mid w\in F\}\cup \{0\}$ and $\overline{u}T^1=\{\overline{uw}\mid w\in F^{1}\}\cup\{0\}$ for $u\in F$. We have two different families of intersections of incomparable principal right ideals for $u,v\in F$:
\begin{align*}
    uT^1\cap vT^1&=(uF^{1}\cap vF^{1})\cup(\{\overline{uw}\mid w\in F\}\cap\{\overline{vw}\mid w\in F\})\cup\{0\};\\
    \overline{u}T^1\cap\overline{v}T^1&=uT^1\cap\overline{v}T^1=(\{\overline{uw}\mid w\in F\}\cap\{\overline{vw}\mid w\in F\})\cup\{0\}.
\end{align*}
By \cite[Corollary 3.6]{Carson}, $F$ is RIH.  Thus, for $u,v\in F$ the intersection $uF^1\cap vF^1$ is generated by some finite set $X$.  It is then straightforward to show that $uT^1\cap vT^1=XT^1$ and $\overline{u}T^1\cap\overline{v}T^1=\{\overline{x} \mid x\in X\}T^1$. Hence $T$ is RIH.

Let $S=T\cup\{a,b\},$
and define a multiplication on $S$, extending that of $T$, by
\[
au=bu=ua=ub=\overline{u}\quad\text{and}\quad at=ta=bt=tb=0\quad (u\in F, t\in\overline{F}\cup\{0\}).
\]
It is straightforward to show that $S$ is a semigroup under this multiplication, and then clearly $S$ is a small extension of $T$. The intersection $aS^1\cap bS^1=\overline{F}\cup\{0\}$ is not finitely generated since $\{\overline{x_i}\mid i\in\mathbb{N}\}$ must be contained in any possible generating set. Hence $S$ is not RIH.

Now we show that the class of RIH semigroups is not closed under large subsemigroups.  Let $M=\{1,g\}\cup S,$ and define a multiplication on $M,$ extending that of $S,$ by
\[
1m=m1=m,\quad g^2=1,\quad ga=ag=b,\quad gb=bg=a,\quad gt=tg=t\quad (m\in M,t\in T).
\]
Then $M$ is a monoid under this multiplication, and clearly $S$ is large in $M.$  We note that the group of units of $M$ is $\{1,g\}\cong\mathbb{Z}_2$.
We claim that $M$ is RIH. Let $u,v\in F\subseteq M$ (all other pairs of elements being $\mathcal{R}$-comparable).

An intersection of incomparable principal right ideals $M$ has one of the following two forms, where $u,v\in F$:
\begin{align*}
    uM\cap vM&=(uF^{1}\cap vF^{1})\cup(\{\overline{uw}\mid w\in F^{1}\}\cap\{\overline{vw}\mid w\in F^{1}\})\cup\{0\},\\
    \overline{u}M\cap \overline{v}M&=\overline{u}M\cap vM=(\{\overline{uw}\mid w\in F^{1}\}\cap\{\overline{vw}\mid w\in F^{1}\})\cup\{0\}. 
\end{align*}
Since $F$ is RIH, for $u,v\in F$ the intersection $uF\cap vF$ is generated by some finite set $X$.  It is then straightforward to show that $uM\cap vM=XM$ and $\overline{u}M\cap\overline{v}M=\{\overline{x} \mid x\in X\}M$. Hence $M$ is RIH.
\end{proof}
\end{counterexample}\vspace{-2mm}

The property of being RIH is preserved under subsemigroups whose complements are ideals, and by small extensions where the complement is an ideal:

\begin{proposition}\label{idealcomplement} Let $S$ be a semigroup with a subsemigroup $T$ such that $S{\setminus}T$ is an ideal of $S$.
\begin{enumerate}[leftmargin=*,topsep=-0.5em,itemsep=-0.5em]
\item If $S$ is RIH then so is $T$.
\item If $S{\setminus}T$ is finite, then $S$ is RIH if and only if $T$ is RIH.
\end{enumerate}
\begin{proof}
(1) Let $a,b\in T$. As $S$ is RIH we have $aS^{1}\cap bS^{1}=XS^{1}$ for some finite set $X$. Since $S{\setminus} T$ is an ideal, it follows that $X\cap T\subseteq aT\cap bT\subseteq (X\cap T)T^1$, and hence $aT\cap bT$ is generated by $X\cap T$.

\noindent (2) Suppose that $T$ is RIH and let $I=S{\setminus}T$. If $a,b\in T$, then $aS\cap bS$ is finitely generated by $X\cup (aI\cap bI)$ where $aT^{1}\cap bT^{1}=XT^{1}$. If $a\in I$ or $b\in I$, then $aS\cap bS\subseteq I$, so $aS\cap bS$ is finite and hence finitely generated.
\end{proof}
\end{proposition}\vspace{-2mm}

\begin{corollary}\label{adjoin0}
A semigroup $S$ is RIH if and only if $S^0$ is RIH.
\end{corollary}\vspace{-2mm}

\begin{counterexample}\label{Howson,ideal}
The class of RIH semigroups is not closed under ideals.
\begin{proof}
The free commutative monoid $\mathbb{N}_0\times\mathbb{N}_0$ on two generators is RIH.  By the proof of \cite[Proposition 4.4]{Carson}, the ideal $\mathbb{N}\times\mathbb{N}$ of $\mathbb{N}_0\times\mathbb{N}_0$ is not RIH.
\end{proof}
\end{counterexample}\vspace{-2mm}

However, if an ideal has an identity element then it forms a retract of the oversemigroup.  Thus, by Corollary \ref{retract}, we have:

\begin{proposition}
\label{retractontoideal} A monoid ideal of an RIH semigroup is RIH.
\end{proposition}\vspace{-2mm} 

It is possible for a semigroup to not be RIH even if both an ideal and the associated Rees quotient {\em are} RIH:

\begin{counterexample} There exists a semigroup $S$ with an ideal $I$ such that $I$ and $S/I$ are RIH but $S$ is not RIH.
\begin{proof}
Let $T$ be any semigroup that is not RIH, and let $\phi:F\rightarrow T$ be a homomorphism from a free semigroup $F$ onto $T$. Let $I=\{a_t\mid t\in T\}\cup \{0\}$ be a set in one-to-one correspondence with $T^{0}$ and let $S=F\cup I$. Define a multiplication on $S$, extending that of $F$, as follows: \[
a_tw=a_{t(w\phi)}\quad\text{and}\quad ab=wa=0w=0\qquad(t\in T,\,w\in F,\,a,b\in I)\] 
It is straightforward to show that $S$ is a semigroup under this multiplication, and clearly $I$ is an ideal of $S$ and also a null semigroup. It it easy to show that $I$ is RIH, and certainly $S/I\cong F^{0}$ is RIH. On the other hand, it is straightforward to show that, for $u,v\in T$, the right ideal $uT\cap vT$ is finitely generated if and only if $a_uS\cap a_vS$ is finitely generated. Since $T$ is not RIH, it follows that $S$ is not RIH.
\end{proof}
\end{counterexample}\vspace{-2mm}

We now turn our attention to direct products.  Recalling the proof of Counterexample \ref{Howson,ideal}, we have:

\begin{counterexample}\cite[Proposition 4.4]{Carson}\label{rih:sgrpdp} The class of RIH semigroups is not closed under direct product.
\end{counterexample}\vspace{-2mm}

On the other hand, the property of being RIH {\em is} preserved under direct factors:

\begin{proposition}\label{productproject}
If the direct product $S\times T$ of two semigroups is RIH, then so are $S$ and $T.$
\begin{proof}
Let $\pi:S\times T\rightarrow S$ be the projection map, and pick any $c\in T.$  Then, for any $\mathcal{R}$-incomparable $a,b\in S,$ we have $(a,c)\pi=a,$ $(b,c)\pi=b$ and
\[\big((a,c)(S\times T)\cap(b,c)(S\times T)\big)\pi=\{(s,t)\mid s\in aS\cap bS, t\in cT\}\pi=aS\cap bS.
\]
Thus, by Proposition \ref{commutes}, $S$ is RIH, and by a dual argument so is $T$.
\end{proof}
\end{proposition}\vspace{-2mm}

Given a semigroup $S,$ we call an element $a\in S$ {\em right factorisable} if $a\in aS$; and we call $S$ {\em right factorisable} if all its elements are right factorisable.  Clearly all monoids and regular semigroups are right factorisable.  By \cite[Proposition 4.5]{Carson}, the direct product of two  right factorisable RIH semigroups is also RIH.  This, together with Proposition \ref{productproject}, yields:

\begin{theorem}\label{rih:prod}
Let $S$ and $T$ be right factorisable semigroups.
Then $S\times T$ is RIH if and only if both $S$ and $T$ are RIH.
\end{theorem}\vspace{-2mm}

We conclude this section by considering free products.

\begin{theorem}\label{rih:sfp}
The semigroup free product of a collection $\{S_i\mid i\in I\}$ of semigroups is RIH if and only if each $S_i$ is RIH.
\begin{proof}
Each $S_i$ is a subsemigroup of $F=\prod^{\ast}\{S_i\mid i\in I\}$ such that $F{\setminus}S_i$ is an ideal of $F$.  Thus, the forward implication follows from Proposition \ref{idealcomplement}(1). The reverse implication is \cite[Proposition 4.2]{Carson}.
\end{proof}
\end{theorem}\vspace{-2mm}

Considering now a monoid free product
$F=\prod^{\ast}_1\{M_i\mid i\in I\}$ of a collection $\{M_i:i\in I\}$ of monoids, \cite[Proposition 4.3]{Carson} states that, if each $M_i$ is RIH, then so is $F$.  The converse also holds, by Proposition \ref{retract}, since each $M_i$ is a retract of $F$; this follows from the well-known universal property that for any monoid $T$ and homomorphisms $\varphi_i : M_i\to T$ ($i\in I$) there exists a homomorphism $\varphi : F\to T$.  Thus, we have:

\begin{theorem}\label{rih:mfp}
The monoid free product 
$\prod^{\ast}_1\{M_i\mid i\in I\}$ of a collection $\{M_i: i\in I\}$ of monoids is RIH if and only if each $M_i$ is RIH.
\end{theorem}\vspace{-2mm}

\section{Finitely right equated semigroups}\label{sec:fre}

\begin{definition}
Let $S$ be a semigroup.  For $a\in S$, the \textit{right annihilator congruence} of $a$ is
\[\bf{r}_{S}(a)=\{(s,t)\in S\times S\mid as=at\}.\]
Clearly $\bf{r}_{S}(a)$ is a right congruence on $S$.  We say that $S$ is \textit{finitely right equated} (\textit{FRE}) if $\bf{r}_{S}(a)$ is finitely generated for each $a\in S$.  
\end{definition}\vspace{-2mm}

\subsection{Basic facts}~

Let $S$ be a semigroup. We denote the identity and universal relations on $S$ by $\Delta_S$ and $\nabla_S$, respectively. Recall that an element $a\in S$ is \textit{left cancellative} if $ab=ac$ implies $b=c$ for all $b,c\in S$, and $S$ is \textit{left cancellative} if each of its elements is left cancellative.

\begin{lemma}\label{delta,nabla} 
For a semigroup $S$ and an element $a\in S$, the following hold.
\begin{enumerate}[leftmargin=*,topsep=-0.5em,itemsep=-0.5em]
\item $\bf{r}_{S}(a)=\Delta_{S}$ if and only if $a$ is left cancellative.
\item $\bf{r}_{S}(a)=\nabla_{S}$ if and only if $aS=\{a^2\}$. Moreover, $\bf{r}_{S}(a)=\nabla_{S}$ and $a$ is right factorisable if and only if $a$ is a left zero.
\end{enumerate}
\begin{proof}
(1) This is clear.
        
(2) If $\bf{r}_{S}(a)=\nabla_{S}$, then for any $s\in S$ we have $(a,s)\in \bf{r}_{S}(a)$ and hence $a^2=as.$ On the other hand, if $aS=\{a^2\}$ then $as=a^2=at$ for all $s,t\in S$, and hence $\bf{r}_{S}(a)=\nabla_S$. The second part of the statement  follows from the first part.
\end{proof}
\end{lemma}\vspace{-2mm}

From Lemma \ref{delta,nabla} we deduce the following corollaries.

\begin{corollary}\label{leftcancellative}
All left cancellative semigroups are FRE.
\end{corollary}\vspace{-2mm}

\begin{corollary}\label{iffuni}
Let $S$ be a semigroup such that, for every element $a\in S$, either $a$ is left cancellative or $aS=\{a^2\}$, and such that there exists at least one $a\in S$ with
$aS=\{a^2\}$. Then $S$ is FRE if and only if $\nabla_{S}$ is finitely generated as a right congruence.
\end{corollary}\vspace{-2mm}

The property that the universal right congruence is finitely generated has recently received significant attention; see for instance \cite{universal, minimal}.

We shall now establish a few useful lemmas.  First, given a relation $\rho$ on a semigroup $S$, we define the \textit{content} of $\rho$ as \[
\bf{C}(\rho)=\{a\in S\mid (a,b)\in \rho{\setminus}\Delta_{S}\text{ for some }b\in S\}.\]
    
\begin{lemma}\label{Cont}
Let $\rho$ be a right congruence on a semigroup $S$. If $\rho=\l X\r$ then the right ideal $\bf{C}(\rho)S^1$ is generated by $\bf{C}(X)$.
\begin{proof}
The statement clearly holds if $S$ is trivial, so assume that $S$ is non-trivial. Certainly $\bf{C}(X)\subseteq\bf{C}(\rho),$ since $X\subseteq\rho$ and so $\bf{C}(X)S^1\subseteq\bf{C}(\rho)S^1$.
Now consider $a\in\bf{C}(\rho)$. Choose $b\in S{\setminus}\{a\}$ such that $(a,b)\in \rho$, so that there exists an $X$-sequence from $a$ to $b.$  It follows that $a\in xS^1$ for some $x\in\bf{C}(X),$ and hence $\bf{C}(X)S^1\supseteq\bf{C}(\rho)S^1$, finishing the proof.
\end{proof} 
\end{lemma}\vspace{-2mm}

\begin{lemma}\label{rightidentity}
Let $S$ be a semigroup.  If there exists a right factorisable element $a\in S$ such that $\bf{r}_{S}(a)$ is finitely generated, then $S=US^1$ for some finite set $U\subseteq S$.
\begin{proof}
By assumption, there exists some $u\in S$ such that $a=au.$
Since $\bf{r}_{S}(a)$ is finitely generated, by Lemma \ref{Cont} the right ideal $\bf{C}(\bf{r}_S(a))S^{1}$ is generated by some finite set $Y.$  Let $U=\{u\}\cup Y.$  We have $(us,s)\in\bf{r}_S(a)$ for all $s\in S.$  Therefore, either $us=s$ or else $s\in\bf{C}(\bf{r}_S(a))S^{1}=YS^1$.  In either case, we have $s\in US^1$.  Thus $S=US^1.$
\end{proof}
\end{lemma}\vspace{-2mm}

\begin{lemma}\label{regelement}
Let $S$ be a semigroup with a regular element $a$, so that there exists $b\in S$ such that $a=aba$.  Then, for any set $U\subseteq S$ such that $S=US^1$ we have $$\bf{r}_{S}(a)=\l\{(bau,u) \mid u\in U\}\r.$$
Consequently, $\bf{r}_{S}(a)$ is finitely generated if and only if $S=US^1$ for some finite set $U\subseteq S$.
\begin{proof}
Let $U\subseteq S$ be such that $S=US^1$, and let $X=\{(bau,u) \mid u\in U\}.$  We claim that $\bf{r}_{S}(a)=\l X\r.$  Indeed, since $abau=au$ for each $u\in U,$ we have $X\subseteq\bf{r}_{S}(a).$  Now let $(s,t)\in\bf{r}_{S}(a).$  By assumption, there exist $u,v\in U$ and $s',t'\in S^1$ such that $s=us'$ and $t=vt'$.  Since $as=at,$ we have an $X$-sequence
$$s=us',\ baus'=bavt',\ vt'=t,$$
as required.  The second part of the lemma follows from the first part and Lemma \ref{rightidentity}.
\end{proof}
\end{lemma}\vspace{-2mm}

\begin{corollary}
\label{regular}
A regular semigroup $S$ is FRE if and only if $S=US^1$ for some finite set $U\subseteq S$.  Consequently, any regular monoid is FRE.
\end{corollary}\vspace{-2mm}

\subsection{Quotients, substructures and extensions}~

Using Corollary \ref{regular}, we may quickly deduce that the property of being FRE is not preserved under removing an adjoined identity.

\begin{counterexample}\label{s1notsR}
If a monoid $S^1$ is FRE then $S$ need not be.
\begin{proof}
Let $S$ be any regular semigroup that is not finitely generated as a right ideal of itself, such as an infinite left zero semigroup.  Then, by Corollary \ref{regular}, the monoid $S^1$ is FRE but $S$ is not FRE.
\end{proof}
\end{counterexample}\vspace{-2mm}

Since left cancellative semigroups are FRE, in particular so are free semigroups. Hence, the class of FRE semigroups is not closed under homomorphic images. Moreover, unlike the case for right ideal Howson semigroups, we do not have closure under Rees quotients.

\begin{counterexample}\label{fre:rq}
The class of FRE semigroups is not closed under Rees quotients.
\begin{proof}
Let $F$ be the free semigroup on an infinite set. Then $F$ is FRE, since it is left cancellative. Let $I$ be the ideal of all words in $F$ of length $2$ or greater.  Then $F/I$ is an infinite null semigroup. From Lemma~\ref{rightidentity}, $F/I$ is not FRE.
\end{proof}
\end{counterexample}\vspace{-2mm}

On the other hand, the property of being FRE is preserved under retracts, as a consequence of the following technical result.  First, given any map $\phi : A\to B$ and any
subset $X\subseteq A\times A,$ we denote by $X\phi$ the set $\{(x\phi,y\phi)\mid(x,y)\in X\}.$

\begin{proposition}\label{homR} 
Let $S$ and $T$ be semigroups such that for each $b\in T$ there exist $a\in S$ and a homomorphism $\phi : S\to T$ such that $a\phi=b$ and $\bf{r}_S(a)\phi=\bf{r}_T(b).$  If $S$ is FRE then so is $T$.
\begin{proof}
Consider $b\in T$. By assumption, there exist $a\in S$ and a homomorphism $\phi : S\to T$ such that $a\phi=b$ and $\bf{r}_S(a)\phi=\bf{r}_T(b).$ 
Since $S$ is FRE, there exists a finite generating set $X$ of $\bf{r}_S(a)$. 
We claim that $X\phi$ generates $\bf{r}_{T}(b).$  Indeed, we have $X\phi\subseteq\bf{r}_S(a)\phi=\bf{r}_{T}(b).$
Now consider $(t,t')\in \bf{r}_{T}(b)$.
Then there exists $(s,s')\in\bf{r}_S(a)$ such that $s\phi=t$ and $s'\phi=t'$. If $s=s'$ then $t=t'$, so assume that $s\neq s'$. Then there exists an $X$-sequence
\[s=p_1c_1,\ q_1c_1=p_2c_2,\ \dots\ ,\ q_nc_n=s'.
\]
By application of $\phi$ we obtain the $X\phi$-sequence
\[t=p_1\phi c_1\phi,\  q_1\phi c_1\phi=p_2\phi c_2\phi,\ \dots\ ,\  q_n\phi c_n\phi=t',\]
interpreting $1\phi$ as $1$.   This establishes the claim, and thus $T$ is FRE.
\end{proof}
\end{proposition}\vspace{-2mm}

\begin{corollary}\label{retractR}
The class of FRE semigroups is closed under retracts.
\begin{proof}
Let $S$ and $T$ be semigroups where $T$ is a retract of $S.$  We show that the condition of Proposition \ref{homR} holds. Let $\phi : S\to T$ be a retraction, and consider $b\in T.$  Then $b\phi=b.$ If $bs=bs'$ then $b(s\phi)=b(s'\phi),$ so $\bf{r}_S(b)\phi\subseteq\bf{r}_T(b).$  Conversely, if $(t,t')\in\bf{r}_T(b)$ then 
clearly $(t,t')\in\bf{r}_S(b),$
and hence $(t,t')=(t,t')\phi \in 
\bf{r}_S(b)\phi$.
\end{proof}
\end{corollary}\vspace{-2mm}

Counterexample \ref{s1notsR} shows that the property of being FRE is not necessarily preserved under large subsemigroups.  We now show 
the same result for monoids.

\begin{counterexample}\label{idealextensionR}\label{fre:lsse}
The class of FRE monoids is not closed under small extensions or large submonoids.
\begin{proof}
Let $F$ be any free commutative monoid on an infinite set. Let $T=\{a\}\cup F^0$, where $a\not\in F,$ and define a multiplication on $T$, extending that of $F^0$, by $a1=1a=a$ and $at=ta=0$ for all $t\in T{\setminus}\{1\}$, where 1 denotes the identity of $F$.  It is straightforward to show that $T$ is a monoid under this multiplication, and obviously $T$ is a small extension of $F.$  The monoid $F$ is cancellative so certainly FRE. On the other hand, since the right ideal $F{\setminus}\{1\}$ of $F$ is not finitely generated, it follows that the right ideal $\bf{C}(\bf{r}_{T}(a))T^1=T{\setminus}\{1\}$ of $T$ is not finitely generated. 
Hence, by Lemma \ref{Cont}, $T$ is not FRE. 

Let $S=T\cup\{e\}$, where $T$ is as above and $e\notin T$, and define a multiplication on $S$, extending that of $T$, by 
$$e1=1e=e,\quad e^2=e,\quad ea=ae=a,\quad eu=ue=0\;\ (u\in F^0).$$
Then $S$ is a monoid under this multiplication, and $T$ is clearly large in $S$.  We claim that $S$ is FRE. Since $e$ is regular, the right annihilator congruence $\bf{r}_{S}(e)$ is finitely generated by Lemma \ref{regelement}. Moreover, we have
\[
\bf{r}_S(a)=(\{1,e\}\times\{1,e\})\cup(S{\setminus}\{1,e\}\times S{\setminus}\{1,e\}).
\]
Since for all $u\in F$ we have $u=1u$ and $0=eu$, it quickly follows that $\bf{r}_S(a)$ is generated by $\{(1,e),(a,0)\}$.
Finally, for all $u\in F^0$ we have
\[
\bf{r}_S(u)=\Delta_S\cup(\{e,a,0\}\times\{e,a,0\}),
\]
which is clearly finitely generated. Hence $S$ is FRE.
\end{proof}
\end{counterexample}\vspace{-2mm}

The following result provides a necessary and sufficient condition for a small extension of an FRE semigroup to be FRE.

\begin{proposition}\label{largeR}
Let $S$ be a semigroup and $T$ a large subsemigroup of $S,$ and suppose that $T$ is FRE.  Then $S$ is FRE if and only if $\bf{r}_S(a)$ is finitely generated for each $a\in S{\setminus}T$.
\begin{proof}
The forward implication is obvious.  For the converse, we only need to show that $\bf{r}_S(a)$ is finitely generated for each $a\in T.$  So, consider $a\in T$.  Let $X$ be a finite generating set for $\bf{r}_{T}(a)$. For each $u\in S{\setminus}T$ such that there exists $v\in T$ with $au=av$, fix some $\alpha_{u}\in T$ such that $au=a\alpha_{u}$. We claim that $\bf{r}_{S}(a)$ is finitely generated by
\[Y=X\cup (\bf{r}_{S}(a)\cap (S{\setminus} T\times S{\setminus} T))\cup \{(u,\alpha_{u})\mid \text{there exists $v\in T$ such that $au=av$}\}.\]
Let $(u,v)\in \bf{r}_{S}(a)$. If $u,v\in T$, then $(u,v)\in\mathbf{r}_T(a),$ so that $u$ and $v$ are connected by an $X$-sequence and hence by a $Y$-sequence. If $u,v\in S{\setminus}T$, then $(u,v)\in Y$. Finally, if, without loss of generality, $u\in S{\setminus}T$ and $v\in T$, then $(u,\alpha_{u})\in Y$ and $(\alpha_{u},v)\in \bf{r}_{T}(a)$. Therefore, there exists an $X$-sequence connecting $\alpha_{u}$ and $v$, and hence a $Y$-sequence connecting $u$ and $v$. This completes the proof.
\end{proof}
\end{proposition}\vspace{-2mm}

Recalling that the right annihilator congruence of the identity of a monoid is trivial, by Proposition \ref{largeR} we have:

\begin{corollary}\label{adjoin1R}
If a semigroup $S$ is FRE then so is $S^1$.
\end{corollary}\vspace{-2mm}

The property of being FRE is inherited by any subsemigroup whose complement is an ideal:

\begin{proposition}\label{idealcomplementR}
Let $S$ be a semigroup with a subsemigroup $T$ such that $S{\setminus} T$ is an ideal of $S$. If $S$ is FRE then so is $T$.
\begin{proof}
Consider $a\in T$, and let $X$ be a finite generating set for $\bf{r}_{S}(a)$. We claim that $\bf{r}_{T}(a)$ is finitely generated by $Y=X\cap(T\times T)$. Let $(u,v)\in \bf{r}_{T}(a)\subseteq \bf{r}_{S}(a)$. Then either $u=v$, in which case we are done, or there exists an $X$-sequence
\[u=p_1c_1,\  q_1c_1=p_2c_2,\ \dots\ ,\  q_nc_n=v.\]
Let $u=q_0$ and $1=c_0$. Consider $0\leq i\leq n-1$, and suppose that $q_i\in T$ and $c_i\in T^1$.  Then $p_{i+1}c_{i+1}=q_ic_i\in T.$  Since $S{\setminus}T$ is an ideal, it follows that $p_{i+1}\in T$ and $c_{i+1}\in T^1$, and $aq_{i+1}=ap_{i+1}\in T,$ so that $q_{i+1}\in T.$  Therefore, by finite induction, we have $p_i,q_i\in T$ and $c_i\in T^{1}$ for all $1\leq i\leq n$. Thus, the above $X$-sequence is a $Y$-sequence, as required.
\end{proof}
\end{proposition}\vspace{-2mm}

\begin{proposition}\label{adjoinzero}
Let $S$ be a semigroup.  Then $S^0$ is FRE if and only if $S$ is FRE and $S=US^1$ for some finite set $U\subseteq S.$
\begin{proof}
Observe that either $S=S^0$ or $S$ is a subsemigroup of $S^0$ whose complement is an ideal.  Thus, using Proposition \ref{idealcomplementR} if necessary, if $S^0$ is FRE then so is $S.$
Now, by Lemma \ref{regelement}, the right annihilator congruence $\bf{r}_{S^{0}}(0)$ is finitely generated if and only if $S^0=V(S^0)^1$ for some finite set $V\subseteq S^0$, which is clearly equivalent to $S=US^1$ for some finite set $U\subseteq S.$  The result now follows from Proposition \ref{largeR}.
\end{proof}
\end{proposition}\vspace{-2mm}

\begin{corollary}\label{adjoin0r}
A monoid $S$ is FRE if and only if $S^0$ is FRE.
\end{corollary}\vspace{-2mm}

\begin{remark}\label{fre:not0}
    It follows from Proposition \ref{adjoinzero} that the class of FRE \textit{semigroups} is not closed under adjoining a zero.
\end{remark}

Moving on to ideals, the proof of Counterexample \ref{s1notsR} yields:
\begin{counterexample}The class of FRE semigroups is not closed under ideals.
\end{counterexample}\vspace{-2mm}

However, if an ideal has an identity element then it forms a retract of the oversemigroup. Thus, by Proposition \ref{retractR}, we have:

\begin{proposition}\label{retractontoidealR} 
A monoid ideal of an FRE semigroup is FRE.
\end{proposition} 
\vspace{-2mm}

A semigroup may not be FRE even if both an ideal and the associated Rees quotient {\em are} FRE:

\begin{counterexample}\label{ISIR}
There exists a semigroup $S$ with ideal $I$ such that $I$ and $S/I$ are FRE but $S$ is not.
\begin{proof} Consider the semigroup 
$T$  as defined in Counterexample \ref{idealextensionR}. Then $I=\{a,0\}$ is an ideal of $T$ that is clearly FRE, and $T/I\cong F^0$ is FRE. However, we have shown $T$ is not FRE.
\end{proof}
\end{counterexample}
\vspace{-2mm}

\subsection{Direct products}~

If a direct product is FRE then so are the direct factors:

\begin{proposition}\label{productprojectR} 
Let $S$ and $T$ be semigroups. If $S\times T$ is FRE then so are $S$ and $T$.
\begin{proof}
Clearly it suffices to prove that $S$ is FRE, and we do so using Proposition \ref{homR}.  Let $\pi : S\times T\to T$ be the projection map, and pick any $b\in T.$  Consider any $a\in S.$  Clearly $\bf{r}_{S\times T}((a,b))\pi\subseteq\bf{r}_S(a),$  and if $(s,s')\in\bf{r}_S(a)$ then $((s,t),(s',t))\in\bf{r}_{S\times T}((a,b))$ for any $t\in T.$  Thus $\bf{r}_{S\times T}((a,b))\pi=\bf{r}_S(a),$ as required. 
\end{proof}
\end{proposition}\vspace{-2mm}

On the other hand, the direct product of two FRE semigroups need not be FRE.

\begin{counterexample}\label{fre:sgrpdp} The class of FRE semigroups is not closed under direct products.
\begin{proof}
Let $S=\{a,0\}$ be the two-element null semigroup and let $T$ any infinite semigroup. Then, for any $x\in T$, any generating set of the right ideal $\bf{C}(\bf{r}_{S\times T}(a,x))S^1=S\times T$ must contain $\{(a,t)\mid t\in T\}$. Therefore, by Lemma \ref{Cont}, $S\times T$ is not FRE.
\end{proof}
\end{counterexample}\vspace{-2mm}

We say that a semigroup $S$ has \textit{pairwise right identities} if for every pair $a,b\in S$ there exists some $s\in S$ such that $as=a$ and $bs=b.$  If $S$ has pairwise right identities then it has {\em finite right identities}, that is, for any finite set $U\subseteq S$ there is an element $s\in S$ such that $us=u$ for all $u\in U$. The following result provides necessary and sufficient conditions for the direct product of two semigroups that have pairwise right identities to be FRE. 

\begin{theorem}\label{productR}
Let $S$ and $T$ be semigroups that have pairwise right identities.  Then $S\times T$ is FRE if and only if $S$ and $T$ are FRE and there exist finite sets $U\subseteq S$ and $V\subseteq T$ such that $S=US^1$ and $T=VT^1$.
\begin{proof}
The direct implication follows from Proposition \ref{productprojectR} and Lemma \ref{rightidentity}.

For the converse, consider $(a,b)\in S\times T$.  By assumption, there exist finite sets $U\subseteq S$ and $V\subseteq T$ such that $S=US^1$ and $T=VT^1$, and $\bf{r}_{S}(a)$ and $\bf{r}_{T}(b)$ are generated by some finite sets $X$ and $Y$, respectively. We claim that $\bf{r}_{S\times T}(a,b)$ is generated by
\begin{align*}
Z=\{((p,p'),(q,q')),\,((u,p'),(u,q')),\,((p,v),(q,v))\hspace{0.5mm}|\hspace{0.5mm}(p,q)\in X,(p',q')\in Y, u\in U,v\in V\}.
\end{align*}
Indeed, let $((s,t),(s',t'))\in \bf{r}_{S\times T}(a,b)$. We may assume that $s\neq s'$.  Since $S$ is FRE, there exists an $X$-sequence
\[s=p_1c_1,\ q_1c_1=p_2c_2,\ \dots\ ,\ q_mc_m=s'.\]
We may assume that $c_i\in S$ for each $1\leq i\leq m$; indeed, as $S$ has pairwise right identities, we may replace any $c_i=1$ with some $d_i\in S$ such that $p_id_i=p_i$ and $q_id_i=q_i$.  
    
Suppose that $t\neq t'$.  Then, since $T$ is FRE and has pairwise right identities, there exists a $Y$-sequence
\[t=p_1'c_1',\ q_1'c_1'=p_2'c_2',\ \dots\ ,\ q_n'c_n'=t'\]
where $c_i'\in T$ for each $1\leq i\leq n$.
Let $N=\max(m,n).$
For each $m\leq i\leq N$ let $q_i=q_m$ and $c_i=c_m$, and for each $n\leq i\leq N$ let $q_i'=q_n'$ and $c_i'=c_n'$.
Then we obtain a $Z$-sequence
$$(s,t)=(p_1,p_1')(c_1,c_1'),\ (q_1,q_1')(c_1,c_1')=(p_2,p_2')(c_2,c_2'),\ \dots\ ,\ (p_N,p_N')(c_N,c_N')=(s',t').$$
Now suppose that $t=t'$.  Since $T=VT^1$ and $T$ has pairwise right identities, it follows that $T=VT$, so that there exist $v\in V$ and $t''\in T$ such that $t=vt''$.  We then have a $Z$-sequence
\[(s,t)=(p_1,v)(c_1,t''),\ (q_1,v)(c_1,t'')=(p_2,v)(c_2,t''),\ \dots\ ,\ (q_m,v)(c_m,t'')=(s',t').\]
This completes the proof.
\end{proof}
\end{theorem}\vspace{-2mm}

\begin{corollary}\label{fre:prod}
Let $S$ and $T$ be monoids.  Then $S\times T$ is FRE if and only if $S$ and $T$ are FRE.
\end{corollary}\vspace{-2mm}

\subsection{Free products}~

We begin this subsection by providing necessary and sufficient conditions for a semigroup free product to be FRE.  

\begin{theorem}
\label{sfp}
Let $\{S_i\mid i\in I\}$ be a collection of semigroups. Then the semigroup free product $F=\prod^{\ast}\{S_i\mid i\in I\}$ is FRE if and only if each $S_i$ ($i\in I$) is FRE and one of the following holds:
\begin{enumerate}[topsep=-0.5em,itemsep=-0.5em]
\item each $S_i$ contains no right factorisable elements; or
\item $I$ is finite, and for each $i\in I$ there exists a finite set $U_i\subseteq S_i$ such that $S_i=U_iS^1$.
\end{enumerate}

\begin{proof}
($\Rightarrow$) For each $i\in I$ the set $F{\setminus}S_i$ is an ideal of $F.$  Thus, by Proposition \ref{idealcomplementR}, each $S_i$ is FRE.

Suppose now that (1) does not hold.  Then there exist some $j\in I$ and $a,u\in S_j$ such that $au=a.$  Then, by Lemma \ref{rightidentity}, there exists a finite set $U_j\subseteq S_j$ such that $S_j=U_jS_j^1$. 
Now, clearly $\{(u\ast\bf{w},\bf{w})\mid\bf{w}\in F\}\subseteq\bf{r}_F(a),$ and it follows that $\bigcup_{i\in I{\setminus}\{j\}}S_i\subseteq\bf{C}(\bf{r}_F(a))$. 
By Lemma \ref{Cont}, the right ideal $\bf{C}(\bf{r}_{F}(a))F^1$ is generated by some finite set $U.$  Since each $F{\setminus}S_i$ is an ideal of $F$, it follows that each $S_i$ ($i\in I{\setminus}\{j\}$) is generated as a right ideal by $U_i:=U\cap S_i$, and hence $I$ is finite.  Thus (2) holds.

$(\Leftarrow)$ Consider an arbitrary element
$\bf{a}=a_1\ast\dots \ast a_n\in F.$  Let $a_n\in S_j$ and let $\rho=\bf{r}_{S_j}(a_n).$  Since $S_j$ is FRE, there exists a finite generating set $X$ of $\rho$.  Let $\rho'$ be the right congruence on $F$ generated by $\rho$ (and hence by $X$), which is easily seen to be   
$$\rho'=\{(s,t),\ (s\ast\bf{w},\ t\ast\bf{w})\mid (s,t)\in\rho,\ \bf{w}\in F\}.$$
If $a_n$ is not right factorisable (in $S_j$) then $\bf{r}_{F}(\bf{a})=\rho'=\l X\r$.  Thus, if (1) holds, then $F$ is FRE.  Suppose now that $a_n$ {\em is} right factorisable.  Then, by assumption, (2) holds.  It is straightforward to show that  
$$\bf{r}_{F}(\bf{a})=
\rho'\cup\{(z\ast \bf{w},\bf{w}),\ (\bf{w},z\ast \bf{w})\mid z\in S_j,\,a_nz=a_n,\ \bf{w}\in F\}.$$ 
Fix $e\in S_j$ such that $a_ne=a_n$.  For each $i\in I$ let $U_i\subseteq S_i$ be a finite set such that $S_i=U_iS_i^1$, and let $U=\bigcup_{i\in I}U_i$.  Note that $U$ is finite (since $I$ is finite) and that $F=UF^1$.
We claim that $\bf{r}_{F}(\bf{a})$ is generated by the finite set
\[Y=X\cup\{(e\ast u,u)\mid u\in U\}.\]
Certainly $Y\subseteq\bf{r}_F(a).$  For the reverse inclusion, since $\rho'=\l X\r$ it suffices to consider a pair $(z\ast \bf{w},\bf{w})$ where $z\in S_j$ is such that $a_nz=a_n$ and $\bf{w}\in F$ is arbitrary.  Clearly $(z,e)\in\rho$, and hence $(z\ast\bf{w},e\ast\bf{w})\in\rho'=\l X\r.$  Now, there exist $u\in U$ and $v\in F^1$ such that $\bf{w}=u\ast v,$ and hence we have a $Y$-sequence $$e\ast\bf{w}=(e\ast u)\ast v,\ u\ast v=\bf{w}.$$
We conclude that there is a $Y$-sequence from $z\ast\bf{w}$ to $\bf{w}$, completing the proof that $\bf{r}_F(\bf{a})=\l Y\r.$  Thus $F$ is FRE.
\end{proof}
\end{theorem}\vspace{-2mm}

\begin{corollary}
Let $\{M_i\mid i\in I\}$ be a collection of monoids.  Then the semigroup free product $F=\prod^{\ast}\{M_i\mid i\in I\}$ is FRE if and only if $I$ is finite and each $M_i$ ($i\in I$) is FRE.
\end{corollary}\vspace{-2mm}

Another immediate consequence of Theorem \ref{sfp} is:

\begin{counterexample}\label{cex:sfp}
Let $S=\{e\}$ be the trivial semigroup, and let $T$ be any semigroup that is not finitely generated as a right ideal of itself.  Then $S\ast T$ is not FRE.
\end{counterexample}\vspace{-2mm}

We now establish that a {\em monoid} free product is FRE if and only if each constituent monoid is FRE.

\begin{theorem}\label{fre:mfp}
Let $\{M_i\mid i\in I\}$ be a collection of monoids.  Then the monoid free product $F=\prod_1^*\{M_i\mid i\in I\}$ is FRE if and only if each $M_i$ ($i\in I$) is FRE.
\begin{proof} 
The direct implication follows from Proposition \ref{retractR} and the fact that each $M_i$ is a retract of $F.$

For the converse, suppose that each  $M_i$ is FRE.  
Consider $\bf{x}=x_n\ast \dots \ast x_1\in F$; we must show that $\bf{r}_F(\bf{x})$ is finitely generated. To this end, define \[N=\min\{i\in I\mid \text{$x_i$ is not right invertible in $M_{x_i\sigma}$}\},\] 
where $\sigma:\bigsqcup_{i\in I}S_i\rightarrow I$ is the source function.  Denote the identity of each $M_i$ by $1_i$, and define a map $$\theta : \bigcup_{i\in I}M_i\to\{1_i\mid i\in I\}$$ by letting $M_i\theta=\{1_i\}$ for each $i\in I.$
For $1\leq i\leq N-1$, let \[R_i=\{r\in M_{x_i\sigma}\mid x_ir=x_i\theta(=1_{x_i\sigma})\}\]
be the set of right inverses of $x_i$ in $M_{x_i\sigma}$. By assumption, for each $1\leq i\leq N$ we may choose a finite set $X_{i}$ of generators for $\bf{r}_{M_{x_i\sigma}}(x_i)$. Fix some $t_i\in R_i$ ($i\in I$) and let
\[X=\bigcup_{i=1}^{N}\{(t_1\ast\dots\ast t_{i-1}\ast p,t_1\ast\dots\ast t_{i-1}\ast q )\mid (p,q)\in X_{i}\}.\]
We claim that $\bf{r}_F(\mathbf{x})=\l X\r$.  Certainly $X\subseteq \mathbf{r}_F(\mathbf{x}),$ and hence $\langle X\rangle\subseteq \mathbf{r}_F(\mathbf{x})$. It remains to prove the reverse inclusion. 
For convenience, we will write $\bf{a}\sim\bf{b}$ if $(\bf{a},\bf{b})\in\langle X\rangle$.

We begin by showing that 
\begin{equation}\label{old1}
t_1\ast t_2\ast\cdots \ast t_\ell\ \sim\ s_1\ast s_2\ast\cdots \ast s_\ell\qquad\text{for any $s_k\in R_u$ for $1\leq k\leq \ell\leq N-1$}.
\end{equation}
The base case, where $\ell=1$, follows immediately from the fact that $(t_1,s_1)\in \mathbf{r}_{M_{x_1\sigma}}(x_1).$ 
Suppose for finite induction that $1\leq \ell<N-1$ and
\[
t_1\ast t_2\ast\cdots \ast t_\ell\,\sim\,s_1\ast s_2\ast\cdots \ast s_\ell.
\]
There is an $X_{\ell+1}$-sequence from
$t_{\ell+1}$ to $s_{\ell+1}$ in $M_{x_{\ell+1}\sigma}$.
Let us write this sequence as 
\[
t_{\ell+1}=p_1c_1,\ q_1c_1=p_2c_2,\ \dots\ ,\ q_hc_h=s_{\ell+1}
\]
where $(p_v,q_v)\in X_{\ell+1}$ and $c_v\in M_{x_{\ell+1}\sigma}$ for $1\leq v\leq h$. Putting
$\mathbf{w}=t_1\ast t_{2}\ast\cdots \ast t_\ell$, it follows that there is an $X$-sequence
\[
\mathbf{w}\ast t_{\ell+1}=(\mathbf{w}\ast p_1)\ast c_1,\ (\mathbf{w}\ast q_1)\ast c_1=(\mathbf{w}\ast p_2) \ast c_2,\ \dots\ ,\ (\mathbf{w}\ast q_h)\ast c_h=\mathbf{w}\ast s_{\ell+1};
\]
that is, 
\[
t_1\ast t_{2}\ast\cdots \ast t_\ell \ast t_{\ell+1}\,\sim\,t_1\ast t_2\ast\cdots\ast t_\ell\ast  s_{\ell+1}.
\] 
An application of the inductive hypothesis and the fact that $\langle X\rangle$ is a right congruence yields 
\[
t_1\ast t_{2}\ast\cdots \ast t_\ell\ast t_{\ell+1}\,\sim\,s_1\ast s_2\ast\cdots\ast s_\ell\ast  s_{\ell+1},
\]
as required.

Next, we note that for any  $1\leq \ell\leq N-1$, if $y\in M_{x_{\ell+1}\sigma}$ and $x_{\ell+1}y=x_{\ell+1}$, then
\begin{equation}\label{S}
t_1\ast t_{2}\ast\cdots \ast t_\ell\ast y\,\sim\,t_1\ast t_2\ast\cdots\ast t_\ell.
\end{equation}
For, given the fact $x_{\ell+1}y=x_{\ell+1}$, we have $(y,x_{\ell+1}\theta)\in \mathbf{r}_{x_{\ell+1}\sigma}(x_{\ell+1})$, and hence there is an $X_{\ell+1}$-sequence in $M_{x_{\ell+1}\sigma}$ from $y$ to $x_{\ell+1}\theta$.
Replacing each instance of $x_{\ell+1}\theta$ with $1$, and then multiplying every term of the resulting sequence on the left by $t_1\ast\cdots \ast t_\ell$, we obtain an $X$-sequence
from $t_1\ast\cdots\ast t_\ell\ast y$ to $t_1\ast\cdots\ast t_\ell$. 

Now, let $\bf{a},\bf{b}\in F$, with reduced forms $\bf{a}=a_1\ast\dots\ast a_p$ and $\bf{b}=b_1\ast\dots\ast b_q$, be such that 
$(\bf{a},\bf{b})\in \mathbf{r}_F(\mathbf{x})$, i.e.\ $\bf{x}\ast\bf{a}=\bf{x}\ast\bf{b}$. We must show that
$(\mathbf{a},\mathbf{b})\in \langle X\rangle$.

The  product $\bf{x}\ast\bf{a}$ reduces to one of the following types of reduced forms:
\[\begin{array}{cll}
\text{(a) }&  x_n\ast\dots\ast x_\ell\ast a_\ell\ast\dots \ast a_p & \text{ if }x_ia_i=1\mbox{ for }1\leq i\leq \ell-1\mbox{ and }x_\ell\sigma\neq a_\ell\sigma;\\
\text{(b) }& x_n\ast\dots\ast x_{\ell+1}\ast x_\ell a_\ell \ast a_{\ell+1}\ast\dots \ast a_p & \text{ if }x_ia_i=1\mbox{ for } 1\leq i\leq \ell-1\mbox{ and }x_\ell a_\ell\neq 1_{x_\ell\sigma}.\end{array}\]
Similarly,  $\bf{x}\ast\bf{b}$ reduces to one of the following types of reduced forms:
\[\begin{array}{cll}
\hspace{-0.8em}\text{(a)}&  \hspace{-0.6em}x_n\ast\dots\ast x_m\ast b_m\ast\dots \ast b_q & \text{if }x_ib_i=1\mbox{ for }1\leq i\leq m-1\mbox{ and }x_m\sigma\neq a_m\sigma;\\
\hspace{-0.5em}\text{(b) }& \hspace{-0.6em}x_n\ast\dots\ast x_{m+1}\ast x_mb_m \ast b_{m+1}\ast\dots \ast b_q & \text{if }x_ib_i=1\mbox{ for } 1\leq i\leq m-1\mbox{ and }x_mb_m\neq 1_{x_m\sigma}.\end{array}\]

We proceed by examining the different cases that can arise.

\noindent{\em Case (i)}.  Suppose that $\bf{x}\ast\bf{a}$ and $\bf{x}\ast\bf{b}$ are each of type (a). Assume without loss of generality that $\ell \geq m$. If $\ell>m$ then we obtain that $a_\ell=x_{\ell-1}$; since $a_{\ell-1}\in R_{\ell-1}$, this contradicts the fact that $\mathbf{a}$ is reduced. Thus
$\ell=m$ and $\mathbf{w}:=a_\ell\ast\dots\ast a_p= b_m\ast\dots \ast b_q$. We then have 
\[\mathbf{a}=a_1\ast\cdots\ast a_{\ell-1}\ast \mathbf{w}\,\sim\,t_1\ast\cdots\ast t_{\ell-1}\ast \mathbf{w}\,\sim\,b_1\ast\cdots\ast b_{\ell-1}\ast \mathbf{w}=\mathbf{b},\]
as required. 

\noindent{\em Case (ii)}. Suppose that  $\bf{x}\ast \bf{a}$ and $\bf{x}\ast \bf{b}$ are each of type (b).  Assume without loss of generality that $\ell \geq m$. This case splits into two subcases.

If $\ell=m$, then we deduce that $(a_\ell,b_\ell)\in \mathbf{r}_{M_{x_\ell}\sigma}(x_{\ell})$ and 
$a_{l+1}\ast\dots \ast a_p=b_{m+1}\ast\dots \ast b_q $. By a now familiar method, from an $X_\ell$-sequence from
$a_\ell$  to $b_\ell$ in $M_{x_\ell\sigma}$ we construct an $X$-sequence from $t_1\ast\cdots\ast t_{\ell-1}\ast a_{\ell}$ to
$t_1\ast\cdots\ast t_{\ell-1}\ast b_{\ell}$ and hence an $X$-sequence from $\mathbf{a}$ to $\mathbf{b}$ in $F$.

Now assume that $\ell>m$, and let $r=l-m-1$. We then have
\begin{equation}\label{new3}
x_\ell a_{\ell}=x_{\ell},\ a_{\ell+k}=x_{\ell-k}\ (1\leq k\leq r),\ a_{\ell+r+1}=x_mb_m,\ 
a_{\ell+r+2}\ast \cdots \ast a_p= b_{m+1}\ast\cdots\ast b_q.
\end{equation}

We will show that
\begin{equation}\label{new4}
a_1\ast\cdots\ast a_{\ell+r}\,\sim\,t_1\ast\cdots\ast t_m.
\end{equation}
To establish (\ref{new4}), we prove by induction that
\[
a_1\ast\cdots\ast a_{\ell+r}\,\sim\, t_{1}\ast\cdots\ast t_{\ell-k-1}\ast w_k\quad\text{for each}\quad0\leq k\leq r,
\]
where $w_k=a_{\ell+k+1}\ast\cdots\ast a_{\ell+r}$ for $0\leq k\leq r-1$ and $w_{r}=1$. Consider first the base case $k=0$. By (\ref{old1}) we have
\[
a_1\ast\cdots\ast a_{\ell+r}\,=\, (a_1\ast\cdots\ast a_{\ell-1})\ast a_{\ell}\ast\cdots\ast a_{\ell+r}\,\sim\, (t_1\ast\cdots\ast t_{\ell-1})\ast a_\ell\ast\cdots\ast a_{\ell+r}.
\]
Since $x_\ell a_\ell=x_\ell$, we 
also have by (\ref{S}) that
\[
a_1\ast\cdots\ast a_{\ell-1}\ast a_{\ell}\,\sim\,t_1\ast\cdots\ast t_{\ell-1}\ast a_{\ell}\,\sim\,t_1\ast\cdots\ast t_{\ell-1}.
 \]
 and so 
\[a_1\ast\cdots\ast a_{\ell+r}\, \sim\, t_1\ast\cdots\ast t_{\ell-1}\ast a_{\ell+1}\ast\cdots\ast a_{\ell+r}.
\]
Now let $0\leq k\leq r-1$ and assume that 
\[
a_{1}\ast\cdots\ast a_{\ell+r}\, \sim\, t_1\ast\cdots\ast t_{\ell-k-1}\ast w_{k}.
\]
Recall that $a_{\ell+k+1}=x_{\ell-k-1}$.  From  $x_{\ell-k-1}t_{\ell-k-1}a_{\ell+k+1}=x_{\ell-k-1}$, we 
may use (\ref{S}) to obtain  \[t_1\ast\cdots\ast t_{\ell-k-2}\ast t_{\ell-k-1}a_{\ell+k+1} \sim  t_1\ast\cdots\ast t_{\ell-k-2}.\] Since $\sim$ is a right congruence, we have
\[
t_1\ast\cdots\ast t_{\ell-k-2}\ast t_{\ell-k-1}\ast w_k\ =\ t_{1}\ast\cdots\ast t_{\ell-k-2}\ast t_{\ell-k-1}a_{\ell+k+1}\ast w_{k+1}\, \sim\, t_1\ast\cdots\ast t_{\ell-k-2}\ast w_{k+1}.
\]
Hence, by the inductive assumption and transitivity, we have $$a_1\ast\cdots\ast a_{\ell+r}\,\sim\,t_1\ast\cdots\ast t_{\ell-(k+1)-1}\ast w_{k+1},$$ completing the inductive step.  Thus (\ref{new4})
is established.

Recalling that $x_mt_mx_m=x_m$, we 
again use (\ref{S}) to obtain  \begin{equation}\label{new5}
t_1\ast\cdots\ast t_{m-1}\ast t_mx_m\,\sim\, t_1\ast\cdots\ast t_{m-1}.  
\end{equation}
We then have
\[
\begin{array}{lllr}
\bf{a}&=&(a_1\ast\cdots\ast a_{\ell+r})\ast x_mb_m\ast b_{m+1}\ast\cdots\ast b_q &\qquad\text{by (\ref{new3})}\\
&\sim& (t_1\ast\cdots\ast t_m)\ast x_mb_m\ast b_{m+1}\ast\cdots\ast b_q &\qquad\text{by (\ref{new4})}\\
&=&(t_1\ast\cdots\ast t_{m-1}\ast t_mx_m)\ast b_m\ast b_{m+1}\ast\cdots\ast b_q &\qquad\\
&\sim& (t_1\ast\cdots\ast t_{m-1})\ast b_m\ast b_{m+1}\ast\cdots \ast b_q&\qquad\text{by (\ref{new5})}\\
&\sim& (b_1\ast\cdots\ast b_{m-1})\ast b_m\ast\cdots\ast b_q&\qquad\text{by (\ref{old1})} \\
&=&\bf{b},& 
\end{array}
\]
as required.

\noindent{\em Case (iii)}. Finally, and without loss of generality, assume that $\bf{x}\ast \bf{a}$ is of type (a) and $\bf{x}\ast \bf{b}$ is of type (b). This case splits into three subcases.

If $\ell> m$, we obtain $a_\ell=x_{\ell-1}b_{\ell-1}$ (if $\ell-1=m$) or $a_\ell=x_{\ell-1}$. Since $a_{\ell-1}\in M_{x_{\ell-1}\sigma}$, this contradicts the fact that $\bf{a}$ is reduced. Thus $\ell\leq m$. If $\ell=m$, then we obtain $x_\ell=x_\ell b_\ell$ and $a_\ell\ast\cdots\ast a_p=b_{\ell+1}\ast\cdots\ast b_q$. From (\ref{S})  we 
have that $t_1\ast\cdots\ast t_{\ell-1} \sim t_1\ast \cdots\ast t_{\ell-1}\ast b_{\ell}$, yielding $\bf{a}\sim \bf{b}$ in $F$. 

Let us consider then the final case, where $\ell< m$.  We have that 
\begin{equation}\label{old5}
x_m=x_mb_m,\quad x_{m-k}=b_{m+k}\ (1\leq k\leq m-\ell )\mbox{ and } a_{\ell}\ast\cdots\ast a_p=b_{2m-\ell +1}\ast\cdots\ast b_q.
\end{equation}
We will show that 
\begin{equation}\label{old6}
t_1\ast\cdots\ast t_{l-1}\,\sim\, t_1\ast\cdots\ast t_{m-1}x_{m-1}\ast x_{m-2}\ast\cdots\ast x_l
\end{equation}
To establish (\ref{old6}), we prove by induction that
$$t_1\ast\cdots\ast t_{l-1}\,\sim\, t_1\ast\cdots\ast t_{l+k}\ast w_k\quad\text{for each}\quad0\leq k\leq m-l-1,$$
where $w_k=x_{l+k}\ast x_{l+k-1}\ast\cdots\ast x_l$. Consider first the base case $k=0.$  Since  $x_{\ell}=x_{\ell}t_\ell x_\ell$
(\ref{S}) gives $t_1\ast\cdots\ast t_{\ell-1} \sim t_1\ast\cdots\ast t_{\ell-1}\ast t_\ell x_\ell=t_1\ast\cdots \ast t_{\ell}\ast w_0$. Now let $0\leq k\leq m-\ell-1$ and assume that $$t_1\ast\cdots\ast t_{\ell-1}\,\sim\, t_1\ast\cdots\ast t_{\ell+k}\ast w_k.$$
Again, since $x_{\ell+k+1}=x_{\ell+k+1}t_{\ell+k+1}x_{\ell+k+1}$, (\ref{S})$ $ gives that 
$$t_1\ast\cdots\ast t_{\ell+k}\ast w_k\, \sim\,  t_1\ast\cdots\ast t_{\ell+k}\ast (t_{\ell+k+1}x_{\ell+k+1})\ast w_k\, =\, t_1\ast\cdots\ast t_{\ell+k+1}\ast w_{k+1}.$$  Hence, by the inductive assumption and transitivity, we have $$t_1\ast\cdots\ast t_{\ell-1}\,\sim\,t_1\ast\cdots\ast t_{\ell+k+1}\ast w_{k+1},$$
completing the inductive step.  Thus (\ref{old6}) is established.

Recalling that $x_m=x_mb_m$, 
by (\ref{S}) we have \begin{equation}\label{old7}
t_1\ast\cdots\ast t_{m-1}\, \sim\, t_1\ast\cdots\ast t_{m-1}\ast b_m.  
\end{equation}
Hence, we have
\[
\begin{array}{lllr}
\bf{a} &=& (a_1\ast\cdots\ast a_{\ell-1})\ast b_{2m-\ell+1}\ast\cdots\ast b_q&\qquad\text{by (\ref{old5})}\\
&\sim& (t_1\ast\cdots\ast t_{\ell-1})\ast b_{2m-\ell +1}\ast\cdots\ast b_q&\qquad\text{by (\ref{old1})}\\
&\sim& (t_1\ast\cdots\ast t_{m-1}x_{m-1}\ast x_{m-2}\ast\cdots\ast x_{\ell})\ast b_{2m-\ell +1}\ast\cdots\ast b_q&\qquad\text{by (\ref{old6})}\\
&=&(t_1\ast\cdots\ast t_{m-1})\ast x_{m-1}\ast x_{m-2}\ast\cdots\ast x_{\ell}\ast b_{2m-\ell+1}\ast\cdots\ast b_q &\\
&\sim& (t_1\ast\cdots\ast t_{m-1}\ast b_m)\ast x_{m-1}\ast\cdots\ast x_{\ell}\ast b_{2m-\ell+1}\ast\cdots\ast b_q&\qquad\text{by (\ref{old7})}\\
&=& (t_1\ast\cdots\ast t_{m-1})\ast b_m\ast b_{m+1}\ast\cdots\ast b_{2m-\ell}\ast b_{2m-\ell +1}\ast\cdots\ast b_q&\qquad\text{by (\ref{old5})} \\
&\sim& (b_1\ast\cdots\ast b_{m-1})\ast b_{m}\ast\cdots\ast b_q&\qquad\text{by (\ref{old1})}\\
&=&\bf{b}.& \\
\end{array}
\]
This completes the proof that $\bf{r}_{F}(\bf{x})=\l X\r.$  Thus $F$ is FRE.
\end{proof}
\end{theorem}\vspace{-2mm}

We now deduce a necessary and sufficient condition for a semigroup free product with an identity adjoined to be FRE.

\begin{corollary}
\label{sfpia}
Let $\{S_i\mid i\in I\}$ be a collection of semigroups, and let $F=\prod^{\ast}\{S_i\mid i\in I\}$.
Then $F^1$ is FRE if and only if each $S_i^1$ ($i\in I$) is FRE.
\begin{proof}
For each $i\in I,$ let $M_i$ be the monoid obtained from $S_i$ by adjoining an element $1_i\notin S_i$ to $S_i$ and defining $1_is=s1_i=s$ for all $s\in S_i$.  Then $F^1\cong\prod_1^{\ast}\{M_i\mid i\in I\}.$  Therefore, by Theorem \ref{fre:mfp}, the monoid $F^1$ is FRE if and only if each $M_i$ is FRE.  Let $i\in I.$  If $S_i$ is not a monoid, then $M_i=S_i^1$.  Suppose now that $S_i$ {\em is} a monoid, i.e.\ $S_i=S_i^1$.  Then, since $S_i$ is an ideal of $M_i$, it follows from Proposition \ref{retractontoidealR} that if $M_i$ is FRE then so is $S_i$.  Conversely, if $S_i$ is FRE then, by Proposition \ref{largeR}, so is $M_i$.  This completes the proof. 
\end{proof}
\end{corollary}\vspace{-2mm}

By Corollaries \ref{adjoin1R} and \ref{sfpia}, we have:

\begin{corollary}\label{cor:sfpia}
Let $\{S_i\mid i\in I\}$ be a collection of semigroups, and let $F=\prod^{\ast}\{S_i\mid i\in I\}$.  If each $S_i$ ($i\in I$) is FRE then so is $F^1$.
\end{corollary}\vspace{-2mm}

It follows from Corollary \ref{sfpia} and Counterexample \ref{s1notsR} that the converse of Corollary \ref{cor:sfpia} does not hold.

\section{Weak right coherency}\label{sec:wrc}

\begin{definition} A semigroup $S$ is \textit{weakly right coherent} (\textit{WRC}) if $S^1$ is weakly right coherent.   
\end{definition} 
\vspace{-2mm}

Since $S^1$ is a free cyclic $S^1$-act, the monoid $S^1$ being WRC is equivalent to every finitely generated right ideal of $S$ being finitely presented as an $S^1$-act.
Moreover, recall that $S^1$ is WRC if and only it is RIH and FRE \cite[Corollary 3.3]{coherentmonoids}, and that $S$ is RIH if and only if $S^1$ is RIH (Proposition \ref{adjoin1}).  Thus, we have:
\begin{theorem}\label{wrc}\label{wrc:adjoin1}
For a semigroup $S$, the following are equivalent:
\begin{enumerate}[topsep=-0.5em,itemsep=-0.5em]
\item $S$ is WRC;
\item $S^1$ is WRC;
\item every finitely generated right ideal of $S$ is finitely presented as a right $S^1$-act;
\item $S^1$ is RIH and FRE;
\item $S$ is RIH and $S^1$ is FRE.
\end{enumerate}
\end{theorem}\vspace{-2mm}

\begin{remark}
There exist WRC semigroups that are not FRE.  An example of such a semigroup is any infinite left zero semigroup $S.$  Indeed, $S$ is easily seen to be RIH, and it was shown in the proof of Counterexample \ref{s1notsR} that $S^1$ is FRE but $S$ is not FRE.    
\end{remark}
\vspace{-2mm}

Given Theorem \ref{wrc}, we may henceforth, without loss of generality, focus on monoids.  The positive results that follow in this section will be proved by juxtaposing the corresponding results from Sections \ref{sec:rih} and \ref{sec:fre} and implicitly applying Theorem \ref{wrc}. 

As free monoids are right coherent \cite[Theorem 1.1]{freemonoid} and hence WRC, certainly the class of WRC monoids is not closed under homomorphic images.  In fact, this class is not closed under Rees quotients:

\begin{counterexample}\label{wrc:rq}
The class of WRC monoids is not closed under Rees quotients.
\begin{proof} Let $S$ be the free semigroup  on an infinite set and let $F=S^1$, so that $F$ is the free monoid on the same set. Certainly $F$ is WRC. With $I$ as in Counterexample \ref{fre:rq}, recall that $F/I$, being isomorphic to the infinite null semigroup, is not finitely generated. Let $a\in S/I$ be a non-zero element. Then $\textbf{C}(\textbf{r}_{F/I}(a))F/I=S/I$ may not be finitely generated and hence, by Lemma \ref{Cont} and Theorem \ref{wrc}, $F/I$ is not WRC.
\end{proof}
\end{counterexample}\vspace{-2mm}

On the other hand, marrying Corollaries \ref{retract} and \ref{retractR} yields:

\begin{theorem}\label{wrc:retract} 
The class of WRC monoids is closed under retracts. 
\end{theorem} 
\vspace{-2mm}

\begin{corollary}
A monoid ideal of a WRC monoid is WRC.
\end{corollary} 
\vspace{-2mm}

\begin{counterexample}\label{wrc:large}\label{wrc:lsse}
The class of WRC monoids is not closed under small extensions or large submonoids.
\begin{proof}
Let $F,$ $T$ and $S$ be as given in Counterexample \ref{idealextensionR}, and recall that $T$ is a small extension of $F$ and a large subsemigroup of $S.$  Since $F$ is a free commutative monoid, it is WRC by \cite[Theorem 4.3]{coherentmonoids}. It was shown that $T$ is not FRE, so that, by Theorem \ref{wrc}, $T$ is not WRC. It was also shown that $S$ {\em is} FRE.  Thus, to show that $S$ is WRC, by Theorem \ref{wrc} it remains to be shown that $S$ is RIH.

Recall that $S=\{a,e\}\cup F^0$ with operation, extending that of $F^0$, given by
\[1s=s1=s\ (s\in S),\quad  eu=ue=au=ua=a^2=0\  (u\in F^0{\setminus}\{1\}),\quad  ea=ae=a,\quad  e^2=e.\]
Since the principal right ideals $0S$, $aS$ and $eS$ are finite, we need only show that intersections of the form 
\[uS\cap vS=uF^0\cap vF^0\ (u,v\in F)\]
are finitely generated. As free monoids are RIH, for $u,v\in F$, the intersection $uF\cap vF$ is finitely generated by some finite set $X$. It is then straightforward to show that $uS\cap vS=XS$, as required.
\end{proof}
\end{counterexample}\vspace{-2mm}

\begin{remark}\label{wrc:ideal}
With $F,$ $T$ and $S$ as in Counterexample \ref{idealextensionR}, observe that $I:=T{\setminus}\{1\}$ is an ideal of $S$ and that $T=I^1$.  Thus, the proof of Counterexample \ref{wrc:large} shows that the class of WRC monoids is not closed under ideals.
\end{remark}\vspace{-2mm}

Now, combining Propositions \ref{idealcomplement} and \ref{idealcomplementR}, we obtain:

\begin{theorem}
Let $S$ be a monoid with a submonoid $T$ such that $S{\setminus} T$ is an ideal of $S$. If $S$ is WRC then so is $T$.
\end{theorem}\vspace{-2mm}

Corollaries \ref{adjoin0} and \ref{adjoin0r} together yield:

\begin{theorem}\label{wrc:0}
A monoid $S$ is WRC if and only if $S^0$ is WRC.
\end{theorem}\vspace{-2mm}

\begin{counterexample}
There exists a monoid $S$ with an ideal $I$ such that $I^1$ and $S/I$ are WRC whilst $S$ is not.
\begin{proof}
Let $S$ and $I$ be as given in Counterexample \ref{ISIR}. Since $I^1$ is finite, it is certainly WRC. As $S/I$ is a free commutative monoid with a zero adjoined, by \cite[Theorem 4.3]{coherentmonoids} and Theorem \ref{wrc:0} it too is WRC. On the other hand, it was shown that $S$ is not FRE, so that, by Theorem \ref{wrc}, $S$ is not WRC.
\end{proof}
\end{counterexample}\vspace{-2mm}

\begin{counterexample}\label{wrc:sgrpdp}
The class of WRC semigroups is not closed under direct product.
\begin{proof}
Free monoids are WRC, as already noted, and hence so are free semigroups by definition.
However, in \cite[Proposition 4.4]{Carson} it was shown that $\mathbb{N}\times\mathbb{N}$ is not RIH. Thus, by Theorem \ref{wrc}, the direct product $\mathbb{N}\times\mathbb{N}$ is not WRC.
\end{proof}
\end{counterexample}\vspace{-2mm}

For monoid direct products, by Theorem \ref{rih:prod} and Corollary \ref{fre:prod} we have:

\begin{theorem}\label{wrc:mdp}
Let $S$ and $T$ be monoids.  Then the direct product $S\times T$ is WRC if and only if both $S$ and $T$ are WRC.
\end{theorem}
\vspace{-2mm}

We round off this section by considering semigroup and monoid free products. By Theorems \ref{rih:sfp} and \ref{sfpia} we have:
\begin{theorem}\label{wrc:sfp}
Let $\{S_i\mid i\in I\}$ be a collection of semigroups, and let $F$ be the semigroup free product $\prod^*\{S_i\mid i\in I\}.$  Then $F^1$ is WRC if and only if each $S_i^1$ ($i\in I$) is WRC.
\end{theorem}\vspace{-2mm}

Finally, from Theorems \ref{rih:mfp} and \ref{fre:mfp} we obtain: 

\begin{theorem}\label{wrc:mfp}
Let $\{M_i\mid i\in I\}$ be a collection of monoids.  Then the monoid free product $F=\prod_1^*\{M_i\mid i\in I\}$ is WRC if and only if each $M_i$ ($i\in I$) is WRC.
\end{theorem} 
\vspace{-2mm}

\section{Axiomatisatability}\label{sec:axiom}

In this final section we provide a survey of connections between axiomatisability of certain classes of monoid acts and each of the properties of being RIH, FRE and WRC.  We begin by establishing the necessary definitions in model theory; we direct the reader to \cite{uniform} for further background.  For basic definitions about acts, the reader should consult \cite{kkm:2000}.

Any class of algebras $\mathcal{A}$ of a given fixed type $\tau$ has an
associated first order language $L^\tau$. If $C$ is a member of $\mathcal{A}$ and $\phi$ is a sentence of 
$L^\tau$ then we say that $C$ is a {\em model} of $\phi$ and write $C\models \phi$ if $\phi$ is a true statement
in $C$. Similarly, if $\Pi$ is a set of sentences of $L^\tau$, then we say that  $C$ is a {\em model} of $\Pi$ and write $C\models \Pi$ if
$C\models \phi$ for every $\phi\in \Pi$. 

\begin{definition}
Let $\mathcal{A}$ be a class of algebras of type $\tau$. A subclass $\mathcal{B}$ of $\mathcal{A}$ is \textit{axiomatisable} if there exists a set of sentences $\Pi$ of $L^\tau$ such that for $C\in\mathcal{A}$, we have $C\in\mathcal{B}$ if and only if $C\vDash \Pi$.
\end{definition}\vspace{-2mm}

Equivalent formulations of the properties of being WRC and being FRE in terms of axiomatisability can be deduced from the work of \cite{modelcompanions}, as follows.

First, recall that the right ideals of a monoid $S$ are precisely the subacts of the right $S$-act $S.$   
A right $S$-act $A$ is said to be \textit{$\alpha$-injective} if for every $S$-homomorphism $\theta:I\rightarrow A$ where $I$ is an $n$-generated right ideal of $S$ with $n<\alpha$, there exists an $S$-homomorphism $\psi:S\rightarrow A$ such that $\psi|_{I}=\theta$. (A right $S$-act is \textit{$n$-generated} if it has a generating set of size $n.$)

\begin{theorem}\cite[Theorem 3]{modelcompanions}\label{alpha}
Let $\alpha\leq\aleph_0$ be a cardinal.  The following are equivalent for a monoid $S$:
\begin{enumerate}[topsep=-0.5em,itemsep=-0.5em]
\item the class of $\alpha$-injective right $S$-acts is axiomatisable;
\item the class of $\alpha$-injective right $S$-acts is closed under ultraproduct;
\item for each $n<\alpha,$ the kernel of every $S$-homomorphism from the $n$-generated free right $S$-act to $S$ is finitely generated.
\end{enumerate}
\end{theorem}\vspace{-2mm}

By Theorem \ref{alpha} and \cite[Corollary 4.2]{modelcompanions}, we have:

\begin{theorem}
The following are equivalent for a monoid $S$:
\begin{enumerate}[topsep=-0.5em,itemsep=-0.5em]
\item the class of $\aleph_{0}$-injective right $S$-acts is axiomatisable;
\item the class of $\aleph_{0}$-injective right $S$-acts is closed under ultraproduct;
\item for each $n<\aleph_0,$ the kernel of every $S$-homomorphism from the $n$-generated free right $S$-act to $S$ is finitely generated;
\item $S$ is WRC.
\end{enumerate}
\end{theorem}\vspace{-2mm}

The right annihilator congruences of a monoid $S$ are precisely the kernels of the $S$-endomorphisms of $S.$ 
Indeed, for any element $a\in S$ there is a natural $S$-endomorphism $\lambda_a : S\to S$ given by $s\lambda_a=as,$ the kernel of which is $\bf{r}_S(a).$ Conversely, the kernel of any $S$-endomorphism $\theta$ of $S$ coincides with $\ker\lambda_{1\theta}=\bf{r}_S(1\theta)$. 
Moreover, it follows from \cite[Theorem 2]{normak} that $\ker\lambda_a$ is finitely generated if and only if the principal right ideal $aS$ of $S$ is finitely presented as a right $S$-act.  Combining these facts with Theorem \ref{alpha} in the case $\alpha=2$, we obtain:

\begin{theorem}
The following are equivalent for a monoid $S$:
\begin{enumerate}[topsep=-0.5em,itemsep=-0.5em]
\item the class of $2$-injective right $S$-acts is axiomatisable;
\item the class of $2$-injective right $S$-acts is closed under ultraproduct;
\item every principal right ideal of $S$ is finitely presented as a right $S$-act;
\item the kernel of every $S$-endomorphism of $S$ is finitely generated;
\item $S$ is FRE.
\end{enumerate}
\end{theorem}\vspace{-2mm}

We now turn our attention to the RIH property. Shaheen \cite{s:2012} showed that this property is equivalent to the axiomatisation of the class of \textit{left} acts satisfying a certain property called Condition (W), which is related to \textit{flatness} \cite{bulmanfleming:1990}.

\begin{definition}\label{defn:W}
Let $S$ be a monoid. 
A left $S$-act $A$ satisfies \textit{Condition (W)} if for any $s,t\in S$ and $a,b\in A$, if $sa=tb$ then there exist $p\in sS\cap tS$ and $c\in A$ such that $sa=tb=pc$.
\end{definition}\vspace{-2mm}

Letting $\mathcal{W}_S$ denote the class of left $S$-acts satisfying Condition (W), we have:

\begin{theorem}\cite[Theorem 5.13]{s:2012}\label{thm:axh} The following are equivalent for a monoid $S$:
\begin{enumerate}[topsep=-0.5em,itemsep=-0.5em]
\item the class $\mathcal{W}_S$  is axiomatisable;
\item the class $\mathcal{W}_S$ is closed under ultraproduct;
\item ultrapowers of $S$ lie in $\mathcal{W}_S$;
\item $S$ is RIH.
\end{enumerate}
\end{theorem}\vspace{-2mm}

This connection between the RIH property and a class of left acts is perhaps not surprising, given the very tight connections in the corresponding situation for rings and modules \cite{chase:1960}. 

Given a monoid $S,$ flatness conditions for a left $S$-act $A$ are defined in terms of the set-valued functor $- \;\! {\otimes _S}\,A:\mathbf{Act}_S\rightarrow\mathbf{Set}$, 
where $\mathbf{Act}_S$ denotes the category of right $S$-acts, preserving certain properties of morphisms. Flatness conditions can be formulated in terms of  complicated sequences
of equalities in acts called \textit{tossings}. However, some forms of flatness can be expressed more simply.  For instance, the strongest notion of flatness, now referred to as \textit{strong flatness}, was shown in \cite[Theorem 5.3]{Stenstrom} to be equivalent to two conditions known as (P) and (E). Strongly flat acts are weakly flat, and an act is weakly flat precisely when it is principally weakly flat and satisfies Condition (W) \cite[Proposition 1.2]{bulmanfleming:1990}.

Strongly flat, weakly flat and principally weakly flat acts all exist within the broader class of torsion-free acts; see \cite[Proposition 1.6]{kilp:1986}.
A left $S$-act $A$ is \textit{torsion-free} if for any $a,b\in A$ and for any right cancellative element $s\in S$ the equality $sa=sb$ implies $a=b$.
Using an example appearing in \cite[Remark III.8.3]{kkm:2000}, we show that Condition (W) is independent of being torsion-free.

\begin{example} \label{ex:w}  
Let $S=(\mathbb{N},\cdot)$ be the natural numbers under multiplication and consider the amalgam $A=\mathbb{N}\coprod^{\mathbb{N}{\setminus}1}\mathbb{N}$. Let us write $A=\{ x,y\}\sqcup\mathbb{N}{\setminus}\{1\}$, where $x$ and $y$ represent the two distinct copies of $1$. Then $A$ is not torsion-free because $2x=2=2y$ but $x\neq y$. However, we claim that $A$ satisfies Condition (W).  Indeed, suppose that  $s,t\in S$ and $a,b\in A$ are such that $sa=tb$. If $sa=tb=x$, then we are forced to have $s=t=1$ and $a=b=x$, so that the claim holds with $p=s$ and $c=x$; similarly if $sa=tb=y$. Otherwise, there exists $c\in A$ such that $sa=tb=pc$ where $p=\text{lcm}(s,t)$. 
\end{example}

\end{document}